\documentclass[12pt]{amsart}
\usepackage{geometry}
\usepackage{amsmath}
\usepackage{amsfonts}
\usepackage{amsthm}
\usepackage{amssymb}
\allowdisplaybreaks[2]

\newtheorem{theorem}{Theorem}

\newtheorem{corollary}[theorem]{Corollary}

\newtheorem{other}{\bf Theorem}

\newenvironment{pf}{\noindent{\emph{Proof.}}}{$\Box$}
\newenvironment{Pf}{\noindent{\emph{Proof of}}}{$\Box$}

\numberwithin{equation}{section}

\newcommand{\D}{\mathbb D}
\newcommand{\C}{\mathbb C}
\newcommand{\hol}{\mathcal Hol}
\newcommand{\ig}{\stackrel{\text{def}}{=}}
\DeclareMathOperator{\og}{O}
\DeclareMathOperator{\op}{o}

\begin{document}

\title[Sequences of zeros and weighted superposition operators]
{Sequences of zeros of analytic function spaces and weighted superposition operators}

\author{Salvador Dom\'{\i}nguez}

\author{Daniel Girela}
\address{An\'{a}lisis Matem\'{a}tico, Facultad de Ciencias, Universidad de M\'{a}laga, 29071 M\'{a}laga, Spain}
\email{sdmolina16@hotmail.com} \email{girela@uma.es}
\thanks{This research is supported in part by a grant from \lq\lq El Ministerio de
Econom\'{\i}a y Competitividad\rq\rq , Spain (MTM2014-52865-P) and
by a grant from la Junta de Andaluc\'{\i}a FQM-210.}

\subjclass[2000]{Primary 30H30, 30H20; Secondary 46E15, 47H99}

\keywords{Weighted superposition operator, sequence of zeros, Bloch
function, Bergman spaces, Weighted Banach spaces of analytic
functions}

\begin{abstract} We use properties of the sequences of zeros of certain spaces of
analytic functions in the unit disc $\mathbb D$ to study the
question of characterizing the weighted superposition operators
which map one of these spaces into another. We also prove that for a
large class of Banach spaces of analytic functions in $\mathbb D$,
$Y$, we have that if the superposition operator $S_\varphi $
associated to the entire function $\varphi $ is a bounded operator
from $X$, a certain Banach space of analytic functions in $\mathbb
D$, into $Y$, then the superposition operator $S_{\varphi ^\prime }$
maps $X$ into $Y$.
\end{abstract}

\maketitle

\section{Introduction}
Let $\D=\{z\in\C:|z|<1\}$ denote the open unit disc in the complex
plane $\mathbb C$ and let $\hol(\mathbb D)$ be the space of all
analytic functions in $\D$ endowed with the topology of uniform
convergence in compact subsets.
\par
Given an entire function $\varphi$, the superposition operator
\[S_\varphi:\hol(\D)\longrightarrow\hol(\D)\]
is defined by $S_\varphi(f)=\varphi\circ f$.
\par More generally, if $\varphi$ is an entire function and $w\in
\hol (\mathbb D)$, the weighted superposition operator
\[S_{\varphi , w}:\hol(\D)\longrightarrow\hol(\D)\] is defined by
$$S_{\varphi , w}(f)(z)\,=\,w(z)\,\varphi \big (f(z)\big ),\quad f\in \hol
(\mathbb D),\quad z\in\mathbb D.$$ In other words, $S_{\varphi ,
w}\,=\,M_w \circ S_\varphi $, where $M_w$ is the multiplication
operator defined by $$M_w(f)(z)=w(z)f(z),\quad f\in \hol (\mathbb
D),\quad z\in\mathbb D.$$ Note that $S_\varphi =S_{\varphi , w}$
with $w(z)=1$, for all $z\in \mathbb D$.
\par
The natural questions in this context are: If $X$ and $Y$ are two
normed (metric) subspaces of $\hol(\D)$, for which pair of functions
$(\varphi , w)$, with $\varphi $ entire and $w\in \hol (\mathbb D)$,
 does the operator $S_{\varphi, w}$ map $X$
into $Y$? When is $S_{\varphi , w}$ a bounded (or continuous)
operator from $X$ to $Y$?
\par Let us remark that we are dealing with non-linear operators.
Consequently, boundedness and continuity are not equivalent a
priori. However, Boyd and Rueda \cite{BR-14} have shown that for a
large class of Banach spaces of analytic functions $X$ and $Y$,
bounded superposition operators from $X$ to $Y$ are continuous.
\par
These questions have been studied for  different pair of spaces $(X,
Y)$, specially in the case of superposition operators. Let us
mention that C\'{a}mera \cite{C} considered the case when $X$ and
$Y$ are Hardy spaces. C\'{a}mera and Gim\'{e}nez \cite{CG}
characterized the entire functions $\varphi $ so that $S_\varphi $
maps a Bergman space into another. The superposition operators
acting between various spaces of Dirichlet type were studied in
\cite{BFV,BuV}. Superposition operators between weighted Banach
spaces of analytic functions were studied in \cite{BV,BR-13,Ramos}.
\'{A}lvarez, M\'{a}rquez, and Vukoti\'c \cite{AMV} studied the
superposition operators between a Bergman space and the Bloch space
in both directions.
\par If $X$ is a subspace of $\hol (\mathbb D)$, a sequence of
points $\{ z_k\} $ in $\mathbb D$ is said to be an \lq\lq
$X$-zero-sequence\rq\rq \, or a \lq\lq sequence of zeros of
$X$\rq\rq \, if there exists a function $f\in X$ which vanishes
precisely in that sequence (counting multiplicities).
\par It has been shown in  \cite{GGM,GM} that information on
the zero-sequences of certain spaces $X, Y$ of analytic functions in
$\mathbb D$ can be an useful tool to characterize the superposition
operators mapping $X$ into $Y$. Our main objective in this paper is
to find new applications of these ideas in the more general context
of weighted superposition operators. We shall present these
applications and some further results in sections~2 and 3.
\par\medskip
\section{Weighted superposition operators between the Bloch space and Bergman spces}
\par For $\,0<p<\infty \,$ and $\,\alpha >-1\,$
the weighted Bergman space $\,A^p_\alpha \,$ consists of those
$\,f\in \hol (\mathbb D)\,$ such that
$$\Vert f\Vert _{A^p_\alpha }\,\ig\, \left ((\alpha +1)\int_\mathbb D(1-\vert z\vert ^2)^{\alpha }\vert f
(z)\vert ^p\,dA(z)\right )^{1/p}\,<\,\infty .$$ Here, $\,dA\,$
stands for the area measure on $\,\mathbb D$, normalized so that the
total area of $\,\mathbb D\,$ is $\,1$. Thus
$\,dA(z)\,=\,\frac{1}{\pi }\,dx\,dy\,=\,\frac{1}{\pi
}\,r\,dr\,d\theta $. The unweighted Bergman space $\,A^p_0\,$ is
simply denoted by $\,A^p$. We refer to \cite{DS,HKZ,Zhu} for the
notation and results about Bergman spaces.
\par A function $f\in \hol (\mathbb D)$ is said to be a Bloch
function if
\begin{equation*}\left \Vert f\right \Vert _{\mathcal
B}\,\ig\,\vert f(0)\vert \,+\,\sup_{z\in \mathbb D}(1-\vert z\vert
^2)\vert f^\prime (z)\vert <\infty .\end{equation*} The space of all
Bloch functions is called the Bloch space and it will be denoted by
$\mathcal B$. We mention \cite{ACP} as a classical reference for
Bloch functions.
\par One of the main results in \cite{AMV} was showing that
if $0<p<\infty $ and $\varphi $ is an entire function, then the
superposition operator $S_\varphi $ maps the Bergman space $A^p$
into the Bloch space $\mathcal B$ if and only if $\varphi $ is a
constant function. The proof of this result given in \cite{AMV} uses
in a very precise way properties of the conformal mapping of the
disc onto a sector. We shall use completely different ideas to prove
the following more general result.

\begin{theorem}\label{ApalphaBloch} Suppose that $0<p<\infty $ and $\alpha >-1$. Let
$w$ be a non-identically zero analytic function in $\mathbb D$ and
let $ \varphi $ be an entire function with $\varphi \not\equiv 0$.
Then the weighted superposition operator $S_{\varphi ,w}$ maps
$A^p_{\alpha }$ into the Bloch space $\mathcal B$ if and only if
$w\in \mathcal B$ and $\varphi $ is constant.
\end{theorem}
\par\medskip When $\alpha =0$ and $w\equiv 1$,
Theorem\,\@\ref{ApalphaBloch} reduces to the just mentioned result
of \cite{AMV}. Before embarking in the proof, we shall recall some
results about the sequences of zeros of the spaces $A^p_\alpha $ and
$\mathcal B$. Horowitz proved in \cite{Ho} the following result.
\begin{other}\label{zeroApalpha} Suppose that $0<p<\infty $ and
$\alpha >-1$. \begin{itemize}\item[(i)] If $g\in A^p_\alpha $,
$g(0)\neq 0$, and $\{ z_k\} $ is its sequence of zeros then
\begin{equation}\label{Ho-zeros-Apa}\prod_{k=1}^n\frac{1}{\vert
z_k\vert }\,=\,\og \left (n^{(1+\alpha )/p}\right ).\end{equation}
\item[(ii)] The statement (i) is sharp in the following sense: For
any given $\varepsilon >0$, there exists a function $g\in A^p_\alpha
$ with $f(0)\neq 0$ whose sequence of zero zeros $\{ z_k\} $
satisfies
\begin{equation}\label{Sharp-Ho-zeros-Apa}\prod_{k=1}^n\frac{1}{\vert
z_k\vert }\,\neq\,\og \left (n^{(1+\alpha )/(p(1+\varepsilon
))}\right ).
\end{equation}
\end{itemize}
\end{other}
\par\medskip Actually, Horowitz proved these results for $\alpha >0$
and Sedletskii \cite{Sed} proved that they remain true for all
$\alpha >-1$. Arguing as in the proof of \cite[Theorem\,\@1]{GNW},
we obtain that $\og \left (n^{(1+\alpha )/p}\right )$ can be
replaced by $\op \left (n^{(1+\alpha )/p}\right )$ in
(\ref{Ho-zeros-Apa}).
\par\medskip Girela, Nowak and Waniurski proved in \cite{GNW} that
if $f$ is a Bloch function with $f(0)\neq 0$ and  $\{ \xi _k\} $ is
its sequence of zeros then
\begin{equation}\label{zerosBloch}\prod_{k=1}^n\frac{1}{\vert
\xi_k\vert }\,=\,\og \left ((\log n)^{1/2}\right ).\end{equation}
Later, Nowak \cite{N} proved that $\og \left ((\log n)^{1/2}\right
)$ cannot be replaced by $\op \left ((\log n)^{1/2}\right )$.
\par\medskip Now we are in disposition to prove
Theorem\,\@\ref{ApalphaBloch}.\par\medskip
\begin{Pf}{\,\em{Theorem\,\@\ref{ApalphaBloch}.}} It is trivial that if
$\varphi $ is constant and $w\in \mathcal B$ then $S_{\varphi
,w}(A^p_\alpha )\,\subset \,\mathcal B$.\par If $\varphi $ is
constant, not identically $0$, and $S_{\varphi ,w}$ maps
$A^p_{\alpha }$ into $\mathcal B$ then it is clear that $w\in
\mathcal B$.
  \par Suppose now that
$w\not\equiv 0$, $\varphi $ is not constant, and $S_{\varphi
,w}(A^p_\alpha )\,\subset \,\mathcal B$. Take $a\in \mathbb C$ such
that $\varphi (a)\neq 0$ and let $f$ be the constant function
defined by $f(z)=a$, for all $z\in \mathbb D$. Since $f\in
A^p_\alpha $, it follows that
$$S_{\varphi ,w}(f)\,=\,\varphi (a)\cdot w\,\in \mathcal B.$$
This implies that $w\in \mathcal B$.
\par Next we take $\varepsilon >0$ and use
Theorem\,\@\ref{Ho-zeros-Apa}\,\@(ii) to pick a function $g\in
A^p_\alpha $ with $g(0)\neq 0$ such that its sequence of zeros $\{
z_k\} $ satisfies (\ref{Sharp-Ho-zeros-Apa}). It is clear that
$\varphi \circ g$ is not constant. Set $F\,=\,S_{\varphi
,w}(g)\,-\,\varphi (0)\cdot w$. Since $w$ and $S_{\varphi ,w}(g)$
are Bloch functions, it follows that
$$F\,=\,S_{\varphi ,w}(g)\,-\,\varphi (0)\cdot w\,=\,w\left [\varphi \circ
g\,-\,\varphi (0)\right ]\,\in \,\mathcal B.$$ Now, it is clear that
$F\not\equiv 0$ and that all the zeros of $g$ are zeros of $F$. In
other words, the sequence $\{ z_k\} $ is contained in the sequence
$\{ \xi _k\} $ of zeros of $F$. Since $\{ z_k\} $ satisfies
(\ref{Sharp-Ho-zeros-Apa}), it follows that that $\{ \xi _k\} $ does
not satisfy (\ref{zerosBloch}). This is in contradiction with the
fact that $F\in \mathcal B$.
\end{Pf}
\par\medskip
Let us turn now our attention to the weighted superposition
operators which map a weighted Bergman space into the Bloch space.
We shall prove the following result which is related to
\cite[Theorem\,\@3]{AMV}.
\begin{theorem}\label{weigtedsuperApalphaBloch} Suppose that:
\begin{itemize}\item [(i)] $0<p<\infty $ and $\alpha >-1$. \item[(ii)] $w\in A^p_{\beta }$ for some
$\beta $ with $-1<\beta <\alpha $. \item[(iii)] $\varphi $ is an
entire function of order less than one, or of order one and type
zero.\end{itemize} Then the weighted superposition operator
$S_{\varphi ,w}$ is a bounded operator from  $\mathcal B$ into
$A^p_{\alpha }$.
\end{theorem}
\begin{pf} Set $M(r)\,=\,\max_{\vert z\vert \le
r}\vert f(z)\vert $ ($0<r<\infty $). Condition (iii) implies that
\begin{equation}\label{logM}\frac{\log M(r)}{r}\,\to \,0,\quad \text{as $r\to \infty
$}\end{equation} (see \cite[Chapter\,\@14]{Hi}). Take $K>0$ and let
$f$ be a Bloch function with norm at most $K$. It is easy to see
(see \cite{ACP}, for example) that there exists an absolute constant
$C>0$ such that
\begin{equation}\label{Bnorm}\vert f(z)\vert \,\le \,\Vert
f\Vert_{\mathcal B}\log\frac{C}{1-\vert z\vert }.\end{equation} By
(\ref{logM}), there exists $r_0>0$ (which depend only on $\varphi $
and $K$) such that
$$\frac{\log M(r)}{r}\,\le \,\frac{\alpha -\beta }{Kp},\quad r\ge
r_0.$$ Thus, using (\ref{Bnorm}), we see that whenever $\vert
f(z)\vert \ge r_0$ we have
\begin{align*}\left \vert S_{\varphi }(f)(z)\right \vert \,= &\,\left
\vert \varphi (f(z))\right \vert \,\le\, \exp \left (\frac{\alpha
-\beta }{Kp}\vert f(z)\vert \right )\\
\le &\,\exp \left (\frac{\alpha -\beta }{p}\log \frac{C}{1-\vert
z\vert }\right )\\ =\,& D (1-\vert z\vert )^{(\beta -\alpha
)/p},\end{align*} with $D =C^{(\alpha -\beta )/p}$. Then it follows
that \begin{equation}\label{intger0}\int_{\vert f(z)\vert \ge
r_0}(1-\vert z\vert )^\alpha \vert w(z)\vert ^p\vert S_\varphi
(f)(z)\vert ^p\,dA(z)\,\le \,D^p\int_ {\mathbb D}(1-\vert z\vert
)^{\beta }\vert w(z)\vert ^p\,dA(z).
\end{equation}
If $\vert f(z) \vert \,\le \,r_0$ then, by the maximum principle,
$\vert S_{\varphi }(f)(z)\vert \le M(r_0)$. Combining this with
(\ref{intger0}), we obtain \begin{align*}&\int_{\mathbb D}(1-\vert
z\vert )^\alpha \vert w(z)\vert ^p\vert S_\varphi (f)(z)\vert
^p\,dA(z)\\ \le \,&D^p\int_ {\mathbb D}(1-\vert z\vert )^{\beta
}\vert w(z)\vert ^p\,dA(z)\,+\,M(r_0)^p\int_ {\mathbb D}(1-\vert
z\vert )^{\alpha }\vert w(z)\vert ^p\,dA(z).\end{align*} This is a
positive number which depend on $\varphi $, $w$, $p$, $\alpha $,
$\beta $, and $K$, but not on $f$. Thus we have shown that
$S_{\varphi ,w} $ is a bounded operator from $\mathcal B$ into
$A_p^\alpha $.
\end{pf}
\par\medskip
\section{Some further results} The arguments used to prove
Theorem\,\@\ref{ApalphaBloch} actually lead to the following result.
\begin{theorem}\label{WSZeros} Let $X$ and $Y$ be two spaces of
analytic functions in $\mathbb D$ satisfying the following
conditions.
\begin{itemize}
\item[(i)] $X$ contains the constant functions.
\item[(ii)] There exists a function $f\in X$ with $f(0)\neq 0$ whose
sequence of zeros $\{ z_k\} $ is not a sequence of zeros of $Y$.
\end{itemize}
Let $w$ be a non-identically zero analytic function in $\mathbb D$
and let $ \varphi $ be an entire function with $\varphi \not\equiv
0$. Then the weighted superposition operator $S_{\varphi ,w}$ maps
$X$ into $Y$ if and only if $\,w\in Y$ and $\varphi $ is constant.
\end{theorem}
\par\medskip
There are a lot of possible choices of the spaces $X$ and $Y$ in
Theorem\,\@\ref{WSZeros} and each one of these choices leads us to a
concrete result. Let us mention just a couple of them.
\par Just as in \cite{BV}, in this paper a weight $\,v\,$ on $\mathbb D$ will be a
positive and continuous function defined on $\mathbb D$ which is
radial, i.\,\@e. $v(z)=v(\vert z\vert )$, for all $z\in \mathbb D$,
and satisfying that $v(r)$ is strictly decreasing in $[0,1)$ and
that $\lim_{r\to 1}v(r)=0$. For such a weight, the weighted Banach
space $H^\infty _v$ is defined by
\begin{equation*}H^\infty _v\,=\,\left \{ f\in \hol (\mathbb D) :
\Vert f\Vert _{v}\,\ig \,\sup_{z\in \mathbb D}v(z)\vert f(z)\vert
<\infty \right\} .\end{equation*} We refer to \cite{BV} for the
origin and relevance of these spaces and for examples of weights on
$\mathbb D$.
\par We recall that a sequence $\{ z_k\} \subset \mathbb D$ is said
to be a Blaschke sequence if $$\sum (1-\vert z_k\vert )\,<\,\infty
.$$ Blaschke sequences are the zero sequences of any of the Hardy
spaces $H^p$ ($0<p\le \infty $) and also of the Nevannlinna class
$N$ (see \cite{Du:Hp}). Actually, using Jensen's formula it follows
that if $f\in \hol (\mathbb D)$ and $f\not\equiv 0$ then its
sequence of zeros $\{ z_k\} $ satisfies the Blaschke condition if
and only if
$$\sup_{0<r<1}\frac{1}{2\pi }\int_0^{2\pi }\log \vert
f(re^{it})\vert \,dt\,<\,\infty .$$ We can state the following
result. \begin{corollary}\label{corHinftyv} Let $Y$ be a space of
analytic functions in $\mathbb D$ such that all the sequences of
zeros of $\,Y$ are Blaschke sequences and let $v$ be a weight on
$\mathbb D$. Let $w$ be a non-identically zero analytic function in
$\mathbb D$ and let $ \varphi $ be an entire function with $\varphi
\not\equiv 0$. Then the weighted superposition operator $S_{\varphi
,w}$ maps $H^\infty _v$ into $Y$ if and only if $\,w\in Y$ and
$\varphi $ is constant.
\end{corollary}
\begin{pf}
By Theorem\,\@\ref{WSZeros}, it suffices to show that there exists a
function $f\in H^\infty _v$, $f\not\equiv 0$, whose sequence of
zeros $\{ z_k\} $ is not a Blaschke sequence. \par Set $\psi
(r)=1/v(r)$ ($r\in [0,1)$). Then $\psi $ is a positive, continuous,
and increasing function in $[0,1)$, and $\lim_{r\to 1}\psi
(r)\to\infty $. Using Lemma\,\@2 in p.\,\@80 of \cite{DS} we can
choose an increasing sequence of positive integers $\{ n_k\}
_{k=1}^\infty $ such that
$$\sum _{k=1}^\infty r^{n_k}\le \frac{1}{2}\log \phi (r),
\quad 0\le r<1.$$ Then the function $f$ constructed in p.\,\@94 of
\cite{DS} with this sequence $\{ n_k\} $ belongs to $H^\infty _v$
and its sequence of zeros is not a Blaschke sequence.
\end{pf}
\par\medskip When $v$ is the weight defined by $v(z)=\left (\log \frac{e}{1-\vert
z\vert }\right )^{-1}$ ($z\in \mathbb D$), the space $H^\infty _v$
is the space of those functions $f\in \hol (\mathbb D)$ such that
\begin{equation}\label{loggrowth}\vert f(z)\vert \,=\,\og \left (\log \frac{1}{1-\vert z\vert
}\right ),\quad \text{as $\vert z\vert \to 1$}.\end{equation} This
is the space which was called $A^0$ in \cite{GNW} and $H^\infty
_{\log }$ in \cite{GPPR} and in \cite{BY}. The space $H^\infty
_{\log }$ and the Bloch space are closely related. It is well known
that any Bloch function $f$ satisfies (\ref{loggrowth}). Hence,
$\mathcal B\subset H^\infty _{\log }$. Since the function
$f(z)=\log\frac{1}{1-z}$ is a Bloch function, the growth condition
(\ref{loggrowth}) is the best one possible for the Bloch space. Bao
and Ye \cite{BY} have identified the space  $H^\infty _{\log }$ as
the dual of a certain Luecking-type subspace of $A^1$ and they have
shown also that the Bloch space is not dense in $H^\infty _{\log }$.
Using Theorem\,\@\ref{WSZeros} and \cite[Theorem\,\@2]{GNW} we
obtain the following result.
\begin{corollary}\label{HinftylogBloch}
Let $w$ be a non-identically zero analytic function in $\mathbb D$
and let $ \varphi $ be an entire function with $\varphi \not\equiv
0$. Then the weighted superposition operator $S_{\varphi ,w}$ maps
$H^\infty _{\log }$ into the Bloch space if and only if $\,w\in
\mathcal B$ and $\varphi $ is constant.
\end{corollary}
\par\medskip Let us recall that the order and type of an entire
function $\varphi $ are the same as those of its derivative (see,
e.\,\@g., \cite[Chapter\,\@ 14]{Hi}). Bearing this in mind and
looking at the results in
\cite{AMV,BV,BR-13,BR-14,BFV,BuV,CG,C,GGM,GM}, we see that for a
good number of pairs of spaces of analytic functions in $\mathbb D$,
$(X, Y)$, we have that if the superposition operator $S_\varphi $
maps $X$ into $Y$ then so does the operator $S_{\varphi ^\prime }$.
Then it is natural to look for general conditions on the pair $(X,
Y)$ which imply this fact. We shall prove a result of this kind when
$Y$ is one of the weighted Banach spaces $H^\infty _v$.

\begin{theorem}\label{Hvvarphivarphiprime}
Let $v$ be weight on $\mathbb D$ and let $(X, \Vert \cdot \Vert )$
be a Banach space of analytic function in $\mathbb D$. Let $\varphi
$ be an entire function. If the superposition operator $S_\varphi $
is a bounded operator from $X$ into $H^\infty _v$, then
$S_{\varphi^\prime }$ maps $X$ into $H^\infty _v$.
\end{theorem}
\begin{pf} Take $f\in X$ and $z\in \mathbb D$.
\par Suppose first that $\vert f(z)\vert \le 1$. Set $A=\sup _{\vert
\xi \vert\le 1}\vert \varphi ^\prime (\xi )\vert $. Then we have
$$\left \vert S_{\varphi ^\prime }(f)(z)\right \vert \,=\,\vert
\varphi^\prime (f(z))\vert \,\le A\,\le\,\frac{Av(0)}{v(z)}.$$
\par Suppose now that $\vert f(z)\vert \ge 1$ and set $R=2\vert
f(z)\vert $. As usual, if $F$ is an entire function and $0\le
r<\infty $, we set
$$M(r, F)\,=\,\sup _{\vert z\vert =r}\vert F(z)\vert .$$
Using Cauchy's integral formula for the derivative and some simple
estimates, we obtain
\begin{align*}\left \vert S_{\varphi^\prime }(f)(z)\right \vert
\,=\,&\left \vert \varphi^\prime (f(z))\right \vert  =\,\left \vert
\frac{1}{2\pi i}\int_{\vert \zeta \vert
=R}\frac{\varphi (\zeta )}{(\zeta -f(z))^2}\,d\zeta \right \vert \\
\le \,& \frac{1}{2\pi }\,2\pi R\,\frac{M(2\vert f(z)\vert ,\varphi
)}{R^2} =\,\frac{M(2\vert f(z)\vert ,\varphi )}{2\vert f(z)\vert }\\
\le \,&\frac{1}{2}\,M(2\vert f(z)\vert ,\varphi ).
\end{align*}
\newline Take now $\theta \in \mathbb R$ such that
$$M(2\vert f(z)\vert , \varphi )\,=\,\left \vert \varphi
\left (2e^{i\theta }f(z)\right )\right \vert $$ and set
$$g(\xi )\,=\,2e^{i\theta }f(\xi ),\quad \xi \in \mathbb D.$$
Then we have \begin{equation}\label{inqsvarphiprime} \left \vert
S_{\varphi^\prime }(f)(z)\right\vert \,\le \frac{1}{2}\,\vert
\varphi (g(z))\vert \,=\,\frac{1}{2}\left \vert S_\varphi
(g)(z)\right \vert .\end{equation} Now, since $f\in X$ we also have
that $g\in X$ and $\Vert g\Vert =2\Vert f\Vert $. Since $S_{\varphi
}$ is a bounded operator from $X$ into $H^\infty _v$, there exists a
positive constant $L$ such that
$$\left \vert S_\varphi (g)(\xi )\right \vert\,\le \,\frac{L\Vert
f\Vert }{v(\xi )},\quad \xi \in \mathbb D.$$ Using this in
(\ref{inqsvarphiprime}), we obtain
$$\left \vert
S_{\varphi^\prime }(f)(z)\right\vert \,\le \,\frac{L\Vert f\Vert
}{2v(z)}.$$
\par Putting both cases together, we obtain
$$\left \vert
S_{\varphi^\prime }(f)(z)\right\vert \,\le \,\frac{C}{v(z)},
\quad\text{for all $z\in \mathbb D$},$$ with $C=\max \left (Av(0),
\frac{L\Vert f\Vert }{2} \right )$. This gives that
$S_{\varphi^\prime }(f)\in H^\infty _v$.
\end{pf}
\par\medskip If $v$ is a weight on $\mathbb D$, we define $DH^\infty
_v$ as follows
\begin{equation*}DH^\infty _v=\left \{ f\in \hol (\mathbb D) : f^\prime
\in H^\infty _v\right \} .\end{equation*} The space $DH^\infty_v$ is
a Banach space with the norm $\Vert \cdot\Vert _{D,v}$ defined by
$$\Vert f\Vert _{D,v}\,=\,\vert f(0)\vert \,+\,\Vert f^\prime \Vert
_{v}.$$
\par\medskip We have the following result.
\begin{theorem}\label{DHvvarphivarphiprime}Let $v$ be a weight on $\mathbb D$ and let $(X, \Vert \cdot \Vert )$
be a Banach space of analytic function in $\mathbb D$. Let $\varphi
$ be an entire function. If the superposition operator $S_\varphi $
is a bounded operator from $X$ into $DH^\infty _v$, then
$S_{\varphi^\prime }$ maps $X$ into $DH^\infty _v$.
\end{theorem}
\par\medskip This result can be proved with arguments very similar
to those used in the proof of Theorem\,\@\ref{Hvvarphivarphiprime},
we omit the details. Notice that when $v(z)\,=\,(1-\vert z\vert )$\,
the space $DH^\infty _v$ reduces to the Bloch space. Hence, as a
particular case of Theorem\,\@\ref{DHvvarphivarphiprime} we obtain.
\begin{corollary}\label{Blochvarphivarphiprime}
Let $v$ be weight on $\mathbb D$ and let $(X, \Vert \cdot \Vert )$
be a Banach space of analytic function in $\mathbb D$. Let $\varphi
$ be an entire function. If the superposition operator $S_\varphi $
is a bounded operator from $X$ into the Bloch space $\mathcal B$,
then $S_{\varphi^\prime }$ maps $X$ into $\mathcal B$.
\end{corollary}

\bibliographystyle{amsplain}

\end{document}